\documentclass{article}

\usepackage{amsmath, amssymb, amsfonts, amsthm, url, graphicx}
\usepackage[all]{xy}

\newcommand{\comment}[1]{}

\newenvironment{equations}{\begin{eqnarray*}}{\end{eqnarray*}}

\newcounter{temp}
\newcommand{\reiterate}[2]{
\setcounter{temp}{\value{mythm}} 
\setcounter{mythm}{#1} #2
\setcounter{mythm}{\value{temp}}
}

\newtheorem{theorem}{Theorem}[section]
\newtheorem{lemma}[theorem]{Lemma}
\newtheorem{proposition}[theorem]{Proposition}
\newtheorem{corollary}[theorem]{Corollary}

\newtheorem{mythm}{Theorem}
\newtheorem{myproposition}[mythm]{Proposition}
\newtheorem{mycorollary}[mythm]{Corollary}

\theoremstyle{definition}
\newtheorem{definition}[theorem]{Definition}
\newtheorem{remark}[theorem]{Remark}

\newcommand{\Ric}{\ensuremath{\, \mathrm{Ric}}}

\newcommand{\sech}{\mbox{ sech }}

\newcommand{\RR}{\mathbb{R}}
\newcommand{\ZZ}{\mathbb{Z}}

\title{An example of how the Ricci Flow can increase topological entropy}
\author{Dan Jane}    
\date{\today}

\begin{document} 

\maketitle

\begin{abstract}  We give a surface for which the Ricci Flow applied to the metric will increase the topological entropy of the geodesic flow.  Specifically, we first adapt the Melnikov method to apply to a Ricci Flow perturbation and then we construct a surface which is closely related to a surface of revolution, but does not quite have rotational symmetry.  This is done by adapting the Liouville metric representation of a surface of revolution.  The final steps of the Melnikov method require numerical integration.  \end{abstract}

\section{Introduction}

Let $M$ be a closed surface with cotangent bundle $\tau: T^*M \to M$.  Let $\omega$ be the canonical symplectic structure and let the Hamiltonian be \begin{displaymath} H(x, p) := \frac{1}{2} \tilde{g}_x(p,p),\end{displaymath} where $\tilde{g}$ is the dual of the Riemannian metric.  The flow $\phi_t$ associated to the Hamiltonian system $(T^*M, \omega, H)$ is called the \emph{geodesic flow} of $(M, g)$ \cite{paternain-1999-}.

An important dynamical invariant for a geodesic flow is the \emph{topological entropy} \cite{brin-2002-}, a measure of how rich the dynamics within the flow is.  On a surface the entropy is non-zero if and only if there exists a transversal homoclinic point, and in this case there will be exponential growth in the lengths of periodic orbits \cite{katok-1980-}.

The \emph{Ricci Flow} of a metric $g$ is a smooth family of metrics $g_t$ satisfying $$\frac{\partial g_{t}}{\partial t} = -2 \Ric_{g_t}, \, \mathrm{with} \, g_{0} = g. $$  Introduced by R. Hamilton \cite{chow-2004-}, the Ricci Flow is in some sense an analogue of the heat equation for the Riemannian metric \cite[p.2]{hamilton-1995-}. Loosely speaking, the Ricci Flow acts to diffuse the curvature across the manifold - homogenising the Riemannian structure.  

The Ricci Flow is particularly well understood on surfaces, where it simplifies the metric towards a constant scalar metric in a controlled fashion, \cite{chow-1991-, hamilton-1988-}.

In \cite{manning-2004-}, Manning used a result of Katok, Knieper and Weiss \cite{katok-1991-} on Anosov flows to show that as the metric on a surface of everywhere negative curvature evolves under the Ricci Flow, the topological entropy of the geodesic flow is strictly decreasing.  

Decreasing entropy corresponds to a simpler dynamical structure, and at a first instance this is what one might always expect as the Riemmanian structure is being simplified.  However, here we give an example of a surface for whom the Ricci Flow will increase the entropy: `simplifying' the metric can break symmetries in the overlying symplectic structure and induce complicated dynamics.

For ease throughout, we call a Riemannian 2-manifold a \emph{surface}, and a compact surface whose geodesic flow is integrable and has a homoclinic connection a \emph{good surface}.  We let the curvature of our good surface be $K$, let $\phi _t$ be the geodesic flow, let $\{x^1, x^2, p_1, p_2\}$ be standard local co-ordinates of the cotangent bundle, let $\Gamma$ be the homoclinic connection, and let $F$ be the second integral of motion.

When a good surface is perturbed through metrics that satisfy the Ricci Flow equation, one can tell whether or not the geodesic flow still has zero entropy by applying the following theorem:

\begin{myproposition} \label{THMdan} We define \begin{displaymath} M(z) := \int^\infty_{-\infty} \frac{\partial K}{\partial x^i} \frac{\partial F}{\partial p_i} (\phi_t(z)) dt.\end{displaymath} If $M(z)$ changes sign then the Ricci Flow will strictly increase the entropy of the flow. \end{myproposition}

The \emph{Normalised Ricci Flow} is a modification of the Ricci Flow that preserves the volume of the manifold; $$\frac{\partial g_t}{\partial t} = - 2 \Ric_{g_t} + \frac{A}{n}g_t, $$ where $A$ is the average curvature and $n = \dim(M)$.  If, in proposition \ref{THMdan}, we were to use the Normalised Ricci Flow instead of the Ricci Flow the conditions for increasing entropy would be the same.  This is unsurprising as the Normalised Ricci Flow is a space-time reparametrization of the Ricci Flow, and the perturbative Melnikov method is the backbone of the proof.

Whilst many examples of good surfaces were known classically - such as ellipsoids and poisson spheres - explicit calculation of the integral in these cases has proved intractable for us.  We show how to construct new examples of good surfaces for which the hypotheses of proposition \ref{THMdan} hold, and go through the calculations in a specific case:

\begin{mythm} \label{THMsurface} If a torus has a Liouville metric \begin{displaymath} ds^2 = \Big(\frac{1}{(1 + 0.125 \cos x)^2} + g(y) \Big) \Big( dx^2 + dy^2 \Big), \end{displaymath}  where $g(y)$ is a smooth $2 \pi$-periodic function strictly greater than $-1/2$ and equal to $y^2(1-y^2)$ on $[-1, 1]$, then the topological entropy of its geodesic flow will be increased by the Ricci Flow.  \end{mythm}

Since the surface above is a torus, the results of Hamilton \cite{hamilton-1988-} imply that as time under the Ricci Flow tends to infinity the metric tends to the flat metric.  Interestingly then, this is a surface for which the metric has zero entropy for $t=0$ and $t \to \infty$, but is not always zero.

Using a branched cover construction, we outline a related example with genus zero.

\begin{mycorollary} There exists a sphere whose geodesic flow has zero entropy initially, but the metric is such that under the Ricci Flow this entropy will increase. \end{mycorollary}

In subsection \ref{SECTmeln} we outline the Melnikov method and then, in \ref{SECTriccimeln}, apply it to a Ricci Flow perturbation to prove proposition \ref{THMdan}.

Section \ref{SECTrfexample} summarises the argument for constructing a surface on which proposition \ref{THMdan} gives a positive result.  The properties of Liouville surfaces are discussed in section \ref{SECTliouvillesurfaces}, as these generalise surfaces of revolution.  Surfaces of revolution can easily be given a homoclinic connection but fail to be useful for another reason, as outlined in remark \ref{RMKrfsor}.  With a well chosen modification we construct a suitable surface in \ref{SECTexplicitcalc}.

The last steps involve a numerical integration and the standard graphs showing convergence are given.

\section{Melnikov's method and the Ricci Flow} \label{SECTrfmelnikov}

\subsection{The Melnikov method} \label{SECTmeln}

Assume the 4-dimensional Hamiltonian system $(M, \omega, H_0)$ has a heteroclinic or homoclinic orbit.  We will consider the effect of perturbing the Hamiltonian by studying the dynamics in the family of conservative systems $(M, \omega, H_\epsilon)$, where $$H_\epsilon = H_0 + \epsilon H_1 + O(\epsilon^2).$$

We assume that the unperturbed system is integrable, with second integral $F$.  We pick a homoclinic orbit and call it $\Gamma$.

\begin{remark}  All arguments also apply in the case of a heteroclinic rather than a homoclinic connection, but for ease of exposition and without significant loss of generality we assume all connections are homoclinic.  \end{remark}

Melnikov's idea was to measure the transversal distance between the stable and unstable manifolds in the perturbed system, using $\Gamma$ as a reference point.  If the stable and unstable manifolds do not coincide then we use theorem 5 of Robinson \cite{robinson-1988-} for analytic Hamiltonians to show that the topological entropy of the system is positive for all $\epsilon > 0$ sufficiently small.  Alternatively, if we do not want to insist in our system being analytic in the relevant region we could also use the results by Burns and Weiss \cite{burns-1998-}.  In this exposition the dynamics in the regions on which we apply the Melnikov method are analytic.  This method is very well understood: see Robinson \cite{robinson-1988-} and Wiggins \cite{wiggins-1988-}, which include discussions on why the integrals below converge.

Let $Y$ be the first order change in the Hamiltonian vector field.  Precisely, $Y$ is such that if we let $X_Z$ be the vector field associated to a Hamiltonian $Z$ then $$X_{H_\epsilon} = X_{H_0} + \epsilon Y + O(\epsilon^2).$$

\begin{definition} The \emph{Melnikov integral} is \begin{displaymath} M(z) := \int_{-\infty}^{\infty} (YF)(\phi_t(z)) dt. \end{displaymath} \end{definition}

\begin{theorem}[Melnikov \cite{melnikov-1963-}, Robinson \cite{robinson-1988-}] \label{THMmelnikovintegral}  If the Melnikov integral $M(z)$ changes sign, then for all $\epsilon > 0$ sufficiently small the entropy of $(M, \omega, H_\epsilon)$ is positive.  \end{theorem}

\subsection{Proof of proposition \ref{THMdan}} \label{SECTriccimeln}

We let $(M, g_t)$ be a Riemannian two-manifold solution of the Ricci Flow.  Since the Gaussian curvature satisfies $K = 2 \Ric$ (also a function of $t$), $$\frac{\partial g_t}{\partial t} = - K g_t.$$

We refer to the standard literature for results: in particular that the Ricci Flow of a surface exists for $t \in [0,T]$, with $T > 0$.  For more information, including the recent proof of the Geometrization Conjecture, see Chow and Knopf \cite{chow-2004-} and Morgan and Tian \cite{morgan-2006-}.

We consider the family of geodesic flows of $(M, g_t)$ given by the Hamiltonian systems $(T^*M, \omega, \tilde{g}_t(p,p))$ for $t \in [0, T]$, $T > 0$.

Assume the geodesic flow at $t=0$ is integrable with second integral $F$ and that the flow has a homoclinic orbit.  We look to apply Melnikov's method to see whether the initial Ricci Flow breaks the homoclinic orbit and hence the perturbed system has positive entropy.

Because the metric is being flowed on the tangent bundle but we need the norm of covectors the following lemma is necessary.

\begin{lemma} $$\frac{\partial \tilde{g_t}}{\partial t} = K \tilde{g_t}.$$ \end{lemma}

\begin{proof} We know $g_t\tilde{g_t}$ is the identity tensor.  On differentiating and applying the product rule the result is immediate.  \end{proof}

\begin{remark}  The Hamiltonian is a function with domain $T^*M$, and as we will see the Gaussian curvature arises naturally when we consider first order perturbations.  However, the Gaussian curvature has domain $M$; therefore we need to use $\hat{K} := K \circ \tau$, where $\tau$ is the footpoint projection mentioned in the introduction.  \end{remark}

\begin{lemma} For a Ricci Flow perturbation of the geodesic flow, $Y = X_{\hat{K}H_0}$. \end{lemma}

\begin{proof}  \begin{equations} H_\epsilon (x, p) = \frac{1}{2} \tilde{g}_{\epsilon}(p,p) &=& \frac{1}{2} \tilde{g}_0(p,p) + \frac{\epsilon}{2} \frac{\partial}{\partial t} \Big|_{t=0} \tilde{g}_t (p,p) + O(\epsilon^2) \\ &=& H_0(p) \hspace{5mm} + \epsilon \hat{K} H_0(p) \hspace{9.8mm} + O(\epsilon^2). \end{equations} Therefore the perturbed Hamiltonian vector field takes the form $X_\epsilon = X_{H_0} + \epsilon X_{\hat{K} H_0} + O(\epsilon^2)$, which implies that $Y = X_{\hat{K} H_0}$. \end{proof}

We use this to simplify $YF$ in theorem \ref{THMmelnikovintegral}.  Using the poisson brackets induced by $\omega$ \begin{equations} YF &=& X_{\hat{K} H_0}F = \{F, \hat{K} H_0\} \\ &=& \hat{K} \{F, H_0 \} + H_0 \{F, \hat{K} \}\\ &=& H_0 \{F, \hat{K} \} \end{equations} since $F$ is an integral.

\begin{proposition} In local co-ordinates $\{x^1, x^2, p_1, p_2 \}$ on $T^*M$, \begin{displaymath} YF = \frac{\partial K}{\partial x^j} \frac{\partial F}{\partial p_j}. \end{displaymath} \end{proposition}

\begin{proof} Since we can consider the flow on the level set $H_0 = 1$ \begin{equations} H_0 \{F, \hat{K} \} &=& \{F, \hat{K} \} \\ &=& -d\hat{K}(X_F) \\ &=& - dK \circ d\tau (X_F) \\ &=& dK(\varphi ),\end{equations} where $\varphi = -d\tau (X_F)$.  The proposition will be proved if we can find $\varphi$ in local co-ordinates.  We have \begin{displaymath} dF = \frac{\partial F}{\partial x^i} dx^i + \frac{\partial F}{\partial p_j} dp_j, \end{displaymath} but as these co-ordinates are Darboux, \begin{displaymath} X_F = \frac{\partial F}{\partial x^i} \frac{\partial}{\partial p_i} - \frac{\partial F}{\partial p_j} \frac{\partial}{\partial x^j}. \end{displaymath}

The result follows on applying $d\tau$, i.e. by setting the $\frac{\partial}{\partial p_i}$ components to zero.  \end{proof}

Together with the result of Robinson \cite{robinson-1988-}, on how the Melnikov integral implies positive entropy, we have proved proposition \ref{THMdan}:

\reiterate{0}{\begin{myproposition} We define \begin{displaymath} M(z) := \int^\infty_{-\infty} \frac{\partial K}{\partial x^i} \frac{\partial F}{\partial p_i} (\phi_t(z)) dt.\end{displaymath} If $M(z)$ changes sign then the Ricci Flow will strictly increase the entropy of the flow. \end{myproposition}}

\section{Construction of a tractable good surface} \label{SECTrfexample}

\begin{definition} A \emph{Liouville surface} is a surface whose metric is locally of the form $(f(x) + g(y))(dx^2 + dy^2)$ for some $f, g$ real functions with $f+g$ positive. \end{definition}

All the good surfaces mentioned in the introduction were quadratically integrable.  It is a classical result (see Darboux \cite{darboux-1894-}) that such surfaces are Liouville surfaces.  The example surface of this paper is also a Liouville surface.

\subsection{Simplifications on a Liouville surface} \label{SECTliouvillesurfaces}

Throughout we let $A = f(x) + g(y)$, so that $ds^2 = A(dx^2 + dy^2).$  Our standard reference is A. V. Bolsinov and A. T. Fomenko, \cite{bolsinov-2004-}.

The Liouville metric is an isothermal parametrization.  Hence the Gaussian curvature is \begin{equation} \label{EQNk} K = -A^{-1} \Delta \log \sqrt{A} = \frac{f'^2 + g'^2}{2A^3} - \frac{f'' + g''}{2A^2}, \end{equation} and the partial derivatives are \begin{equation} \label{EQNcurvatures} \left\{ \begin{matrix} K_x = \left( -f'''A^2 + 2f'g''A +4f'f''A - 3f'g'^2 - 3f'^3 \right)/(2A^4), \\ K_y = \left( -g'''A^2 + 2f''g'A +4g'g''A - 3f'^2g' - 3g'^3 \right)/(2A^4). \end{matrix} \right. \end{equation}

\begin{lemma}  \label{THMliouvilleint} The geodesic flow on a Liouville surface is integrable, with integral $F = p_x^2 - fH$. \end{lemma}

\begin{proof} By direct calculation. \end{proof}

\begin{lemma} \label{THMorbitequation} For any geodesic $\gamma = (x, y)$ on a Liouville surface, there exists a constant $a$ such that \begin{equation} \label{EQNorbits} \frac{dy}{dx} = \pm \sqrt{\frac{g(y) - a}{f(x) + a}}. \end{equation} \end{lemma}

\begin{proof}  We pass to the tangent bundle, using $p_x = A\dot{x}$, $p_y = A\dot{y}$.  Since $H$ and $F$ are preserved by the flow, so is $$a := \frac{F}{H} = \frac{g p_x^2 - f p_y^2}{p_x^2+p_y^2} = \frac{g \dot{x}^2 - f \dot{y}^2}{\dot{x}^2 + \dot{y}^2}. $$  Now we solve for $\dot{y}/\dot{x}$. \end{proof}

To simplify the exposition we throw away some generality at this stage.  With theorem \ref{THMsurface} in mind, we assume the Liouville metric is defined on a cylindrical neighbourhood arising from a lifted co-ordinate strip $\RR \times [-a,a]$.  

A geodesic on the surface lifts to a line in the co-ordinate strip.  Hence a point $z(t)$ moving along a geodesic on the surface determines a point $(x(t), y(t)) \in \RR \times [-a,a]$.  Once we set one point in the co-ordinate chart, $x(t)$ and $y(t)$ are uniquely defined.  If we assume $\dot{y}/\dot{x}$ is bounded, then each value of $x$ uniquely determines a value of $y$, and we can equivalently describe the geodesic by the function $y = y(x)$.

We make sure that geodesics under consideration are trapped in this neighbourhood, and that the derivative bounded, in remark \ref{RMKtrappy} below.

\begin{corollary} \label{THMLiouIntegral} The Melnikov integral in proposition \ref{THMdan} simplifies on such a Liouville strip; \begin{equation} \int^{\infty}_{-\infty} \frac{\partial K}{\partial x^i} \frac{\partial F}{\partial p^i} (\phi_t(z) ) dt = 2 \int^\infty_{-\infty} K_x(x, y) g(y) - K_y(x, y) f(x) \frac{dy}{dx} \,\, dx. \end{equation}  \end{corollary}

\begin{proof} We substitute $F = p_x^2 - fH$ as a second integral. \begin{equations} \int^{\infty}_{-\infty} \frac{\partial K}{\partial x^i} \frac{\partial F}{\partial p^i} (\phi_t(z) ) dt &=& \int^{\infty}_{-\infty} 2 \frac{K_x g(y)p_x - K_y f(x)p_y}{f+g} (\phi_t(z) ) dt \\ &=& 2 \int^{\infty}_{-\infty} K_x g(y)\dot{x} (\phi_t(z) ) dt - 2 \int^{\infty}_{-\infty} K_y f(x) \dot{y} (\phi_t(z) ) dt. \end{equations} Then we make the substitution $x = x(t)$ and view $y$ as a function of $x$ as we run along the orbit $\phi_t(z)$.  The improper integrals must all converge absolutely for the argument to hold, but this has been assumed since the beginning.  \end{proof}

For clarity we define $I(x,y) := K_x(x, y) g(y) - K_y(x, y) f(x) dy/dx$.

\begin{remark} From now on we restrict the flow to $H^{-1}(1)$, a 3-manifold.  Obviously local co-ordinates are $\{x, y, \theta\}$, where $\tan \theta = p_x / p_y$. \end{remark}

\begin{proposition} \label{THMhypcoordinates} In addition to our current hypotheses, if $g'(0) = 0$ and \begin{displaymath} g''(0) > f'^2(0)/f(0) - f''(0) \end{displaymath} then $y = 0$ is a hyperbolic geodesic. \end{proposition}

\begin{proof}  Liouville's theorem implies the phase space will foliate into 2-manifolds wherever $dF$ and $dH$ are linearly independent.  Direct calculation shows that they are nowhere parallel; the foliation fails only because $dF=0$.  

This foliation cannot hold true on $y=0$ if it is a hyperbolic geodesic, and so the proof will follow on considering conditions for $dF = 0$, which are \begin{displaymath} -f' \sin^2 \theta dx + g' \cos^2 \theta dy - (f + g) \sin 2 \theta d\theta = 0 \end{displaymath}  The last coefficient, $(f + g) \sin 2 \theta$, is only zero when $\sin 2 \theta$ is zero, as $f + g$ is always positive.  Without loss of generality this implies that $\theta = \pi/2$.  The second term yields $g'(0) = 0$.

As $\theta = \pi/2$ at every point along $y=0$, the axis is a geodesic.  By substituting the condition for $g''(0)$ into (\ref{EQNk}) we see that the curvature along $y=0$ is everywhere negative.  Hyperbolicity of the geodesic along the axis then follows from a standard Jacobi field argument.  \end{proof}

\begin{remark} \label{RMKrfsor}  It is well known that the Ricci Flow preserves the isometry group of the metric, see Hamilton \cite[p.3]{chow-2004-}.  Thus if we were to Ricci Flow a surface of revolution, the $S^1$ symmetry would be preserved and the surface would remain a surface of revolution, and so remain integrable with second integral the Clairaut integral.

This implies the geodesic flow on a surface of revolution will remain integrable along the Ricci Flow; on surface flows this implies that the entropy remains zero.  Hence  if there was a homoclinic connection in the geodesic flow, no homoclinic tangle would form: the Melnikov integral on a surface of revolution should always be zero.

A surface of revolution is exactly a Liouville surface with $f \equiv 0$.  By the above argument, for a surface of revolution the hypotheses of proposition \ref{THMdan} should not apply and the Melnikov integral of corollary \ref{THMLiouIntegral} should always be zero:

\begin{equations} M(z) &=& 2 \int^\infty_{-\infty}K_x(x, y) g(y) - K_y(x, y) f(x) \frac{dy}{dx} \, dx \\ &=& 2 \int^\infty_{-\infty} 0 \cdot g(y) - K_y(x, y) \cdot 0 \cdot \frac{dy}{dx} \, dx  = 0, \end{equations} as expected.

In the next subsection we look to break the $S^1$ symmetry by modifying $f \equiv 0$. \end{remark}

\subsection{Explicit calculation of the Melnikov integral} \label{SECTexplicitcalc}

On a good surface, we have a hyperbolic orbit $c$ with a homoclinic connection $\Gamma$.  The geodesic flow is integrable with second integral $F$.

To calculate an explicit local relation $y = \gamma(x)$ which describes a geodesic in $\Gamma$ we must solve the orbit equation in lemma \ref{THMorbitequation}.  By continuity, the value of $a$ on the homoclinic connection is the same as on $c$, the hyperbolic orbit.  We will assume that on $c$ we have $y = 0$ and $p_y = 0$, as in proposition \ref{THMhypcoordinates}. Since we are on the unit tangent bundle these assumptions imply that $a = g(0)$.  

We set $g(0) = 0$ and from (\ref{EQNorbits}) obtain the following equation: \begin{equation} \label{EQNorbit}  \gamma'(x)= \pm \sqrt{\frac{g(y)}{f(x)}}. \end{equation}  The sign determines whether we are considering the trajectory above or below $y=0$; the geodesics $y = \pm \gamma(x)$.  Due to the symmetry of the example outlined below we consider the positive solution without loss of generality.

\begin{figure}[!hbtp]
\begin{center}
\scalebox{0.7}
{
\includegraphics[width=15cm]{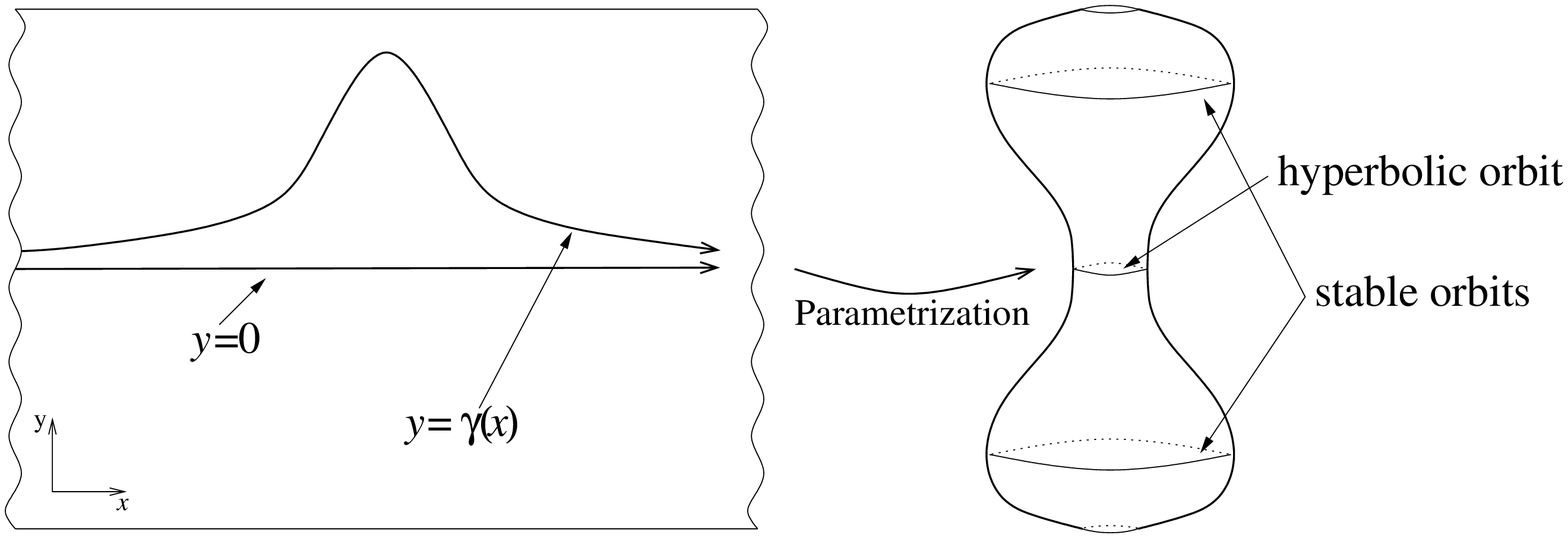}
}
\caption{Parametrization of the surface by a horizontal strip of $\RR \times [-1, 1]$. The local co-ordinates must be periodic in $x$ in order to be able to quotient to a cylinder.}
\label{PICsurface}
\end{center}
\end{figure}

We look to model our surface on a surface of revolution, as in remark \ref{RMKrfsor}, and wish to cover a cylindrical open set of a surface with Liouville co-ordinates.  We will need $f$ periodic in order to take a quotient of $\RR \times [-1, 1]$ and form a cylinder (see figure \ref{PICsurface}).

To ensure $f$ is always positive after periodic modification, we take our Liouville surface with $f \equiv 1$.

\begin{remark}  \label{RMKtrappy} If $f$ is constant, then the metric, function $\gamma$, curvature, and $g$ are all dependent: in particular, we can set $\gamma$ rather than set $g$.  Hence we can ensure that the geodesic is restricted to the co-ordinate chart, $\RR \times [-1, 1]$, and that a value of $x$ uniquely determines a value of $y$ along a geodesic in $\Gamma$.  \end{remark}

We use this observation to construct the surface, and then perturb $f$ away from $1$.  Although various functions are suitable candidates for $\gamma$, we chose $\gamma(x) = 1/2 \sech x$ as the convergences are exponential and the surface is smooth.  This induces $g(y) = y^2(1-y^2)$.  As expected, this is a function with a single minimum - exactly as in the surface of revolution case, remark \ref{RMKrfsor}.

We choose our perturbation of $f$ by asking what properties we require in $y = \gamma(x)$.  From the form of $g$ determined by $\gamma$ in remark \ref{RMKtrappy}, and from (\ref{EQNorbit}), the orbit equation, \begin{displaymath} g(y) = \left( \frac{d}{dy} \gamma^{-1}(y) \right)^{-2} \quad \Rightarrow \quad y = \gamma \left(\kappa + \int \frac{dx}{\sqrt{f(x)}} \right), \end{displaymath} for some constant of integration $\kappa$. The function $f$ must be smooth and periodic - and hence bounded.  The simplest candidate is $f(x) = (1 + \mu \cos x)^{-2}$, for some small $\mu$, so that a geodesic in the homoclinic connection $\Gamma$ is determined by \begin{equation} y = \gamma( \kappa + x + \mu \sin x) = \frac{1}{2} \sech (\kappa + x + \mu \sin x). \end{equation}

\begin{remark} The Melnikov function is evaluated at $z = (x,y)$, a point on the homoclinic orbit.  Since $\kappa$ fixes $z$ and vice versa we consider $M(\kappa)$ instead of $M(z)$.  In the particular example just outlined $M(\kappa)$ is a an odd function of $\kappa$:  the aim is simply to find a value of $\kappa$ such that $M(\kappa)$ is non-zero.  \end{remark}

Unfortunately, attempts to show that integrals of this form are non-zero through analytic means alone have so far been unsuccessful.  We proceed with a method that requires numerical analysis for one part of the argument.

We split the integral into two parts, $M(\kappa) = A_L(\kappa) + B_L(\kappa)$, where \begin{displaymath} A_L(\kappa) = \int_{-L}^L I(x, \gamma(\kappa + x + \mu \sin x)) dx \end{displaymath} and $B_L(\kappa)$ is the integral over the remaining outlying regions.

The argument proceeds as follows.  We fix $\kappa$ and $L$.  We show that $A_L$ is large to within some accuracy $e$, whereas $|B_L|$ is small.  Then $M(\kappa) \geq A_L - e - B_L > 0$, so that $M(\kappa)$ is non-zero and we are done.

\begin{lemma} \label{THMbbound} For $0 < \mu < 1/4$, \begin{displaymath} B_L(\kappa) \leq 1200 (1 + 2 e^{2|\kappa|}) \gamma(L)^2. \end{displaymath} \end{lemma}

\begin{proof} \begin{displaymath} I(x,y) = K_x g - K_y f dy/dx = f dy/dx (K_x dy/dx - K_y). \end{displaymath}  $f$ is bounded by $(1-\mu)^{-2}$ and $dy/dx$ goes to zero exponentially, so it remains to bound the change in curvature.

As an example we bound $K_x$ along $\gamma$; $K_y$ follows in a very similar manner.  We bound $f$, $g$ and their derivatives before plugging them into (\ref{EQNcurvatures}) for $K_x$ and simplifying.

$|f(x)| \leq (1-\mu)^{-2}$, as mentioned, and $|f'(x)| \leq 2 \mu (1 - \mu)^{-3}$ is also trivial.

\begin{displaymath} f''(x) = -2 \mu \frac{2\mu \cos^2 x - \cos x - 3 \mu}{(1 + \mu \cos x)^4}. \end{displaymath}  Given $ 0 < \mu < 1/4$, it is simple to show that the norm of the numerator has maximum $1 + \mu$.  This gives $|f''(x)| \leq 2\mu (1 + \mu) (1 - \mu)^{-4}$.  Analogously $|f'''(x)| \leq 2\mu (1 + 4\mu) (1 - \mu)^{-5}$.

Since $0 < \mu < 1/4$, we have \begin{displaymath} \begin{matrix} |f(x)| \leq 2, & \qquad  & |f'(x)| \leq 2, \\ |f''(x)| \leq 2, & \qquad & |f'''(x)| \leq 5. \end{matrix} \end{displaymath}

Consider the half $x > 0$.  We know that here $y$ is monotonically decreasing on $x$.  The derivatives of $g(y) = y^2(1-y^2)$ can be bounded by plugging in the maximum value of $y$, which is less than $\gamma(L - |\kappa| - \mu) = \gamma_0$.  This $\gamma_0$ is less than 1, as $\gamma(x) \leq 1$ for all $x$, and we use this below.  The argument for $x < 0$ is symmetric.

Thus \begin{displaymath} \begin{matrix} |g(y)| \leq \gamma_0^{\phantom{0}2}, & \qquad  & |g'(y)| \leq 2 \gamma_0, \\ |g''(y)| \leq 2, & \qquad & |g'''(y)| \leq 24 \gamma_0. \end{matrix} \end{displaymath}

We note $A^{-1} = (f(x) + g(y))^{-1} \leq (1 + \mu)^2$.  On substitution into (\ref{EQNcurvatures}), \begin{equations} |K_x| & \leq & 24 \gamma_0 (1+ \mu)^4 + (4 + 8) (1 + \mu)^6 + (6 \gamma_0^{\phantom{0}2} +  4 \gamma_0^{\phantom{0}3})(1 + \mu)^8 \\ & \leq & (1 + \mu)^8 (24 + 12 + 10) \leq 276. \end{equations}  The calculation of $|K_y| \leq 288 \gamma_0$, which tends to zero as $|L + \kappa + 1| \to \infty$,  reflects $K_y(x, 0) = 0$.

Lastly, $|dy/dx| \leq 5\gamma(L)/4$ as $\tanh$ is bounded by 1, and another easy calculation gives $\gamma_0 \leq \gamma(L) \cdot \exp(2\mu + 2 |\kappa|^2)$.

We are now in a position to show \begin{equations} \left| \int_L^\infty I(x,y) dx \right| & \leq & \int_L^\infty |f(x)| \left( |K_x| \cdot \left| \frac{dy}{dx} \right| + |K_y| \right) \left| \frac{dy}{dx} \right| dx \\ & \leq & 2 \left( 276 \gamma(L) + 288 \gamma_0 \right) \int_{\gamma(L)}^{\gamma(\infty)} -dy \\ & \leq & 600 (1 + 2 e^{2|\kappa|}) \gamma(L)^2. \end{equations}

$B_L$ is the sum of two integrals both bounded by this value, which completes the analysis.  \end{proof} 

\begin{remark}  This bound is very rough - perhaps to the casual observer recklessly so - but it won't matter as $\gamma(x) \to 0$ exponentially fast.  For example, with $|\kappa| \leq 3$ and $0 < \mu < 1/4$, we only need $L = 10$ for $B_{10}(\kappa) \leq 0.005$. \end{remark}

\begin{figure}[!hbtp]
\begin{center}
\scalebox{0.7}
{
\includegraphics[width=10cm]{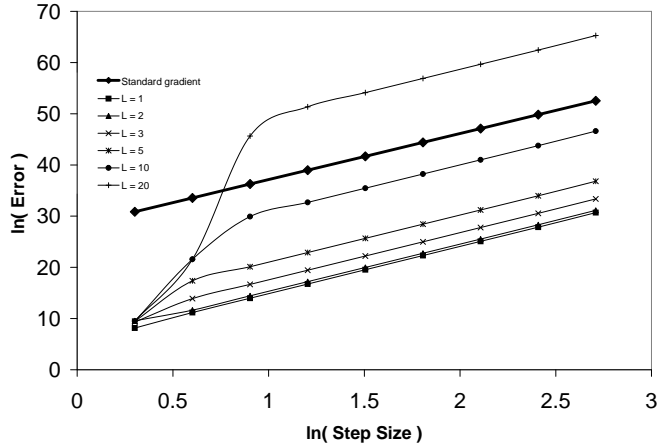}
}
\caption{log$_e$-log$_e$ graph of convergence as step size decreases, with $\kappa = 1$ and $\mu = 1/8$.  Standard convergence is given as a reference gradient.  Initial noise is due to inaccurately low step sizes.}
\label{PICerrors}
\end{center}
\end{figure}

We use Maple to implement the fourth order classical Runge-Kutta method on $A_L(\kappa)$, with $\kappa = 1$ and (arbitrarily) $\mu = 1/8 \neq 0$.  Evidence that the error is being controlled as usual is shown in figure \ref{PICerrors}, which shows convergence as expected.

\begin{figure}[!hbtp]
\begin{center}
\scalebox{0.7}
{
\includegraphics[width=10cm]{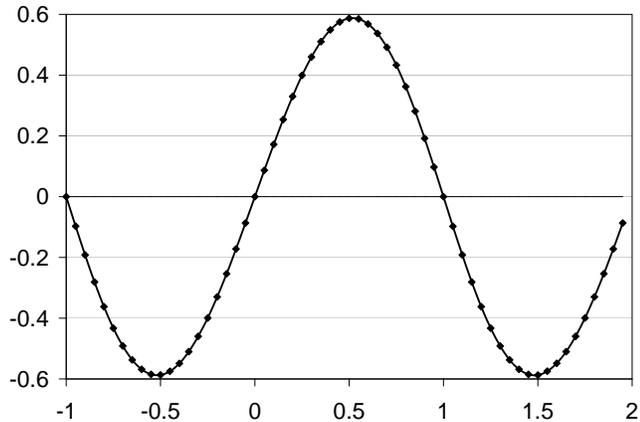}
}
\caption{Melnikov Integral as a function of $\kappa$ ($x$-axis in multiples of $\pi$), with $\mu = 1/8$.}
\label{PICsinusoidal}
\end{center}
\end{figure}

\begin{remark} As expected, if we hold $L$ constant and calculate $A_L(\kappa)$ in terms of varying $\kappa$ (with a reasonable accuracy) we get a sinusoidal graph, figure \ref{PICsinusoidal}.  \end{remark}

Table \ref{TABintegration} shows how the values of $A_L(1)$ converge as $L$ increases, and we also include the bound for $B_L(1)$ calculated in lemma \ref{THMbbound}.  This gives irrefutable evidence that the Melnikov integral is non-zero: $0.479$ to three significant figures.  By proposition \ref{THMdan}, we have shown

\begin{table}
\caption{\label{TABintegration}Values of $A_L(1)$ and the bound for $B_L(1)$ against $L$, to 8 significant figures, with $\mu = 1/8$.}
\begin{center} \begin{tabular}{cr@{.}lr@{.}l} 
\ $L$ & \multicolumn{2}{c}{$A_L(1)$} & \multicolumn{2}{c}{bound for $B_L(1)$} \\ 
\hline 
1 & 1&5990894 & 7&951683E4 \\
2 & 1&0073283 & 1&337684E4 \\
3 & 0&59562584 & 1&868009E1 \\
5 & 0&48440447 & 3&438049 \\
10 & 0&47929061 & 1&561013E-4 \\
20 & 0&47929047 & 3&217488E-13 \\
\hline \end{tabular} \end{center}
\end{table}

\reiterate{1}{\begin{mythm}  If a torus has a Liouville metric \begin{displaymath} ds^2 = \Big(\frac{1}{(1 + 0.125 \cos x)^2} + g(y) \Big) \Big( dx^2 + dy^2 \Big), \end{displaymath}  where $g(y)$ is a smooth $2 \pi$-periodic function strictly greater than $-1/2$ and equal to $y^2(1-y^2)$ on $[-1, 1]$, then the topological entropy of its geodesic flow will be increased by the Ricci Flow.  \end{mythm}}

\begin{proof}[Proof of corollary C]  Quotienting by the involution $\sigma(x,y) = (-x, -y)$ gives a branched double cover of $T^2$ over $S^2$.  We assume $F$ and $G$ are both even functions, so they are invariant under $\sigma$.  If $F$ and $G$ are chosen such that $F(x)+G(y)$ vanishes at the four branch points in a suitable manner, then the induced symmetric tensor on $S^2$ will be a Liouville metric.

The `suitable manner' is that near the zeroes \begin{equation} F(t) + G(it) = 0, \end{equation} and that at any zero $F$ has a non-zero second derivative.  These relations arise since the branch points are double branch points, as in \cite{bolsinov-2004-}.

We take $F$ similar to $f$ above, but subtracting a constant so that it has zeroes; $F(x) = (1+ \mu \cos x)^{-2} - (1 + \mu)^{-2}$.  By adding the same constant to $G$ the metric is unaffected and the argument for the Ricci Flow to increase the entropy still holds.  We shift the strip on which the argument was applied away from the branch points using $\tilde{y} = y - 3$.   Hence a suitable function for $G$ is a smooth function reflecting at $0$ and $3$ such that \begin{displaymath} G(y) = (1+ \mu)^{-2} + \left\{ \begin{matrix}  \tilde{y}^2( 1-\tilde{y}^2) & y \in [2, 3] \\ - (1 + \mu \cosh y)^{-2} & y \in [0, 1] \end{matrix} \right. \end{displaymath}  See figure \ref{PICsphere}.

\begin{figure}[!hbtp]
\begin{center}
\scalebox{0.5}
{
\includegraphics[width=14cm]{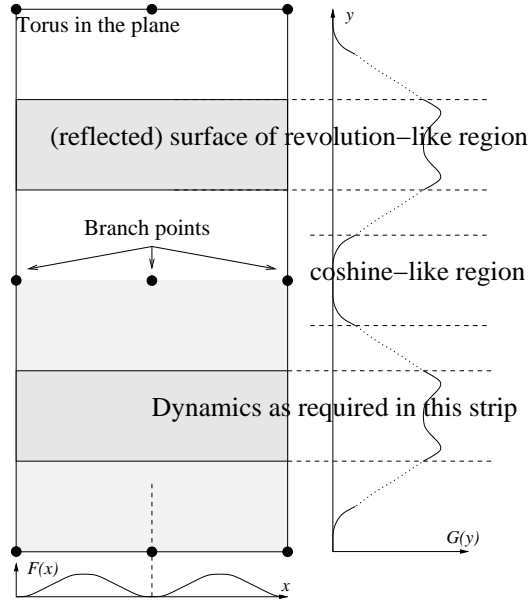}
}
\caption{How $F$ and $G$ are adapted to the branched cover of $T^2$ over $S^2$.  Zeroes of $F$ and $G$ occur simultaneously at the branch points.}
\label{PICsphere}
\end{center}
\end{figure}

Thus we have defined a Riemannian metric on $\RR^2/(2\pi\ZZ \oplus 12\ZZ)/\sigma$, which is topologically a sphere.  On this surface the geodesic flow has zero topological entropy, but under the Ricci Flow this entropy will become strictly positive.  \end{proof}

\section*{Acknowledgements}

I wish to thank my Ph.D. supervisor, Gabriel Paternain, for many helpful discussions, and Leo Butler who suggested investigating Liouville surfaces and ways to improve the exposition. 

I would also like to thank Jonathan Dawes for advice on using and explaining numerical methods, and Martin Kerin, Leonardo Macarini and Anthony Manning for their useful suggestions.

\bibliography{Master}
\bibliographystyle{plain}

\end{document}